\documentclass[a4paper,11pt]{article}

\usepackage{amsthm} 
\usepackage[tbtags]{amsmath} 

\usepackage{amsfonts}
\usepackage{amssymb}
\usepackage{verbatim}

\def\vf{\vspace {4 mm}}

\def\E{\mathbb{E}}

\def\eps{\epsilon}

\def\1{\mathbf{1}}

\def\lam {\lambda}
\def\tce{t_c + \eps}
\def\tce2{t_c + \frac{\eps}{2}}

\def\var{\text{var}}
\def\cov{\text{cov}}

\newtheorem{thm}{Theorem}
\newtheorem{lem}{Lemma}

\newtheorem{prop}{Proposition}

\title{The Forgetfulness of Balls and Bins}

\author{William Perkins \footnote {Courant Institute of Mathematical Science, New York University, New York, NY. Supported in part by NSF grant OISE-0730136 }}

\begin{document}

\maketitle

\begin{abstract}
We find the asymptotic total variation distance between two distributions on configurations of $m$ balls in $n$ labeled bins: in the first, each ball is placed in a bin uniformly at random; in the second, $k$ balls are planted in an arbitrary but fixed arrangement and the remaining $m-k$ balls placed uniformly at random.
\end{abstract}

\section{Introduction}

Planted distributions arise in several contexts in the study of random structures and algorithms.  They are used as a means of studying complicated conditional distributions: instead of studying a random structure conditioned on the presence of a substructure, we can instead plant a copy of the substructure and add random elements around it.  One example is the random graph: we can plant a subgraph, say a $k$-clique or a Hamiltonian cycle, in an empty graph on $n$ vertices, then add random edges on top (see \cite{alon1998finding}, \cite{broder1994finding},\cite{kuera1995expected}).  Another class of examples are random satisfiability problems, including Random $k$-SAT and random graph $k$-coloring. In the case of Random $k$-SAT, instead of conditioning on there being a satisfying assignment to a random formula, we can pick a random assignment uniformly, then sample random clauses from the set of clauses satisfied by that assignment.  Planted Random $k$-SAT is thus guaranteed to have satisfying assignment and has been used to test algorithms and in the analysis of the satisfiability threshold (see for example, \cite{achlioptas2008algorithmic}, \cite{coja2007almost}).

In each of these models, a natural question is how close the planted distribution is to the conditional distribution or to the basic random distribution.   The problem of statistically distinguishing a planted Hamiltonian cycle in a random graph was proposed by Klas Markstr\"{o}m and considered by a group (including Svante Janson, Colin McDiarmid, Oliver Riordan and Joel Spencer) at the Discrete Probability program (Spring 2009) at Institut Mittag-Leffler.  The original problem was to determine how many random edges are needed to``hide" a Hamiltonian cycle on $n$ vertices by adding $m-n$ random edges.  The planted Hamiltonian cycle is hidden if the planted distribution and standard random graph with $m$ edges are statistically indistinguishable asymptotically, i.e. the total variation distance between the two distributions tends to $0$ as $n \to \infty$.  

In this paper we consider a `pure' version of the problem: a planted version of the balls-and-bins model. We determine the number of random balls needed to `hide' an initial planted configuration of balls. We have two primary motivations for studying this model: first, we can answer the question of the distinguishability of the planted distribution completely, finding the distinguishability threshold and total variation limit for any given planting; second, many discrete probability problems can be reduced to an instance of balls-and-bins, and so a balls-and-bins result may prove to be a useful tool.

The standard balls-and-bins model involves throwing $m$ balls into $n$ bins, with each ball independently thrown into a uniformly chosen bin.  The standard model induces a probability distribution on configurations of $m$ balls in $n$ bins which we will call the $ST(ANDARD)$ distribution.  In the planted version of the model, we begin with a fixed arrangement of balls already in bins: perhaps one ball planted in each bin, or $k$ balls planted in the first bin and none in any others.  If we have planted $k$ balls, we then throw the remaining $m-k$ balls into bins uniformly at random.  The planted model induces its own distribution on configurations of $m$ balls in $n$ bins which we will call the $PL(ANTED)$ distribution.  The $PL$ distribution depends on the particular initial planting.  Our main question is: how large must $m$ be as a function of $n$ so that the total variation distance between $ST$ and $PL$, $||ST -PL||_{TV}$, tends to $0$?  In other words, how many random balls do we need to throw in bins to ``forget'' our initial planting?

\section{Preliminaries}

Here we introduce our notation for the paper.  A configuration $Z$ of $m$ balls in $n$ bins consists of a list of non-negative integers, $\{z_i\}_{i=1}^n$, $\sum z_i =m$, where $z_i$ is the number of balls in bin $i$.  An initial planting $A$, of $k$ balls in $n$ bins consists of $\{a_i \}_{i=1}^n$, $\sum a_i =k$, where $a_i$ is the number of balls planted in bin $i$.  In the labeled case, these lists are ordered, while in the unlabeled case the lists are unordered.  $PL(Z)$ is the probability of configuration $Z$ under the ST distribution, and $PL_A(Z)$ or $PL(Z)$ is the probability under the planted distribution with initial planting $A$ (with the subscript omitted if the initial planting is clear from the context).  We measure the distinguishability of the two distributions by their total variation distance:
\begin{equation}
\label{tvdef}
 ||ST - PL||_{TV} = \frac{1}{2} \sum_{Z} |ST(Z) - PL(Z) | 
\end{equation}
where the sum is over all possible configurations $Z$ of $m$ balls in $n$ bins. The total variation distance depends on $m$, $n$, and the initial planting $A$.  We write $TV_A(m,n)$ for the total variation distance between $ST$ and $PL$ with $m$ balls, $n$ bins and initial planting $A$, or simply $TV(m,n)$ if the initial planting is understood.

We write \lq with high probability\rq\ or \lq whp\rq\ if an event holds with probability  $\to 1$ as $n \to \infty$ (other authors sometime use `aas' or `asymptotically almost surely'). We use the standard asymptotic notation $O(\cdot), o( \cdot), \Theta(\cdot), \Omega(\cdot),$ and $\omega(\cdot)$.  We often combine the asymptotic notation with a statement about probability: $g(n) = O(f(n))$ whp means that there exists a constant $K$ so that $\Pr[g(n) \le K f(n)] \to 1$ as $n \to \infty$.  $g(n) = o(f(n))$ means that for any constant $c>0$, $\Pr[g(n) > c f(n)] \to 0$ as $n \to \infty$.

\section{Main Result}
Our main result divides the set of initial plantings into three regimes and characterizes the asymptotic behavior of $TV(m,n)$ in each case.  For a given initial planting $A$, we define $V(A)$ as follows: 

\begin{equation}
\label{Vdef}
 V(A) = \frac{ \sum a_i^2 }{n } - \frac{k^2}{n^2} 
\end{equation}

We can view $V(A)$ as the variance of the random variable $\mathcal{A}$ that takes value $a_i$, $i= 1, \dots n$, with probability $\frac{1}{n}$.  It is a characterization of how `flat' or `hilly' the initial planting is.

\begin{thm}
 \label{mainthm}
Let $A$ be an initial planting of $k$ balls in $n$ labeled bins, with $k >> \sqrt n$.   Then:
\begin{enumerate}
\item The Flat regime: $V(A) = o\left ( \frac{k}{n^{3/2}} \right)$.  Let $m \sim c k n^{1/2}$.  Then
\[ TV(m,n) =  2 \Phi \left ( \frac{1}{2 \sqrt{2}c}  \right) -1   +o(1) \]
where $\Phi(\cdot)$ is the standard normal distribution function.
 \item The Hilly regime: $V(A) = \omega\left ( \frac{k}{n^{3/2}} \right)$.  Let $m \sim c  V n^2$, then
\[ TV(m,n) = 2 \Phi \left ( \frac{1}{2 \sqrt{c}}  \right) -1 + o(1)\]
\item The Intermediate regime: $V(A) \sim \lam \frac{k}{n^{3/2}}$ with $\lam$ a constant.  Let $m \sim c k n^{1/2}$.  Then
\[TV(m,n) =  2\Phi \left( \frac{1}{2} \sqrt{\frac{\lam}{c} +\frac{1}{2c^2} }\right) -1 +o(1) \]
\end{enumerate}
\end{thm}

We make a few remarks about the theorem. In each regime, as $c \to \infty$, $TV(m,n) \to 0$ and as $c \to 0$, $TV(m,n) \to 1$, so this scaling gives the correct threshold at which the two distributions become distinguishable.  There is also a smooth transition between the three regimes.  

Note that for some values of $k$, only Hilly arrangements exist.  For example, for $k=o(n)$, all initial arrangements are Hilly.   Also, the Flat and Intermediate regimes are quite restrictive. For $k=n$, for example,  $n - O(\sqrt n)$ of the $a_i$'s are exactly $1$ in both the Flat and Intermediate regimes.  A randomly chosen initial planting will be Hilly with high probability.   

\paragraph{Illustrative Special Cases of the Theorem}

Two special cases of the theorem will be useful guides to what follows.
\begin{enumerate}
 \item The pure flat initial planting, starting with $1$ ball in each bin, or $a_i =1$ for all $i$.  In this case the theorem states that the distinguishability threshold occurs at $m =\Theta(n^{3/2})$
\item A pure hilly planting, starting with $k$ balls in the first bin, $a_1 =k$, $a_i =0$, $i =2, \dots n$.  In this case the theorem states that the distinguishability threshold occurs at $m = \Theta(k^2 n)$.  

\end{enumerate}

The remainder of the paper is devoted to the proof of Theorem \ref{mainthm} and is organized as follows: in Section \ref{lowerboundsec} we prove the lower bounds on total variation distance, and introduce distinguishing statistics in each regime.  In Section \ref{upperboundsec} we prove the upper bounds.  Sections \ref{normalsec}, \ref{errorsec}, and \ref{calcsec} are devoted to technical lemmas, and in Section \ref{concludesec} we give some concluding remarks.

\section{Lower Bounds}
\label{lowerboundsec}
Consider the following game with the goal of distinguishing between the ST and PL distributions: one of the distributions is chosen at random with probability $\frac{1}{2}$ each and then a configuration $Z$ sampled from it.  The player of the game sees only the configuration and must determine which distribution it came from.  He wants a strategy that maximizes the probability of selecting the correct distribution.  For example, in the case of planting all $k$ balls in bin $1$, one natural strategy would be to look at $z_1$, the number of balls that end up in bin $1$.  If $z_1$ is higher than some threshold, choose the PL distribution, otherwise choose the ST distribution.  

Such a strategy gives a lower bound to $TV(m,n)$ in the following way: via Bayes' formula we see that the optimal strategy would be to compute $ST(Z)$ and $PL(Z)$ and choose whichever is higher. If we call the probability of success using the optimal strategy $p^*$ then we have 
\begin{equation*}
  p^* = \frac{1}{2} \sum_{ST(Z) \ge PL(Z)} ST(Z) + \frac{1}{2} \sum_{PL(Z) > ST(Z)} PL(Z)
\end{equation*}
And so
\begin{equation*}
 2 p^* -1 = \sum_{ST(Z) \ge PL(Z)} ST(Z) + \sum_{PL(Z) > ST(Z)} PL(Z) - \sum_Z PL(Z)
\end{equation*}
\begin{equation*}
 = \sum_{ST(Z) \ge PL(Z)} [ST(Z) - PL(Z)] = TV(m,n)
\end{equation*}
So given any strategy with probability of success $p$, we have $p^* \ge p$ and  

\begin{equation}
 \label{tvlowerboundeq}
TV(m,n) \ge 2p -1
\end{equation}

In this section we will give strategies in each of the three regimes, calculate their success probabilities and find the lower bounds for Theorem \ref{mainthm}.  The strategies are all similar to the strategy described above for the extreme case: we choose a statistic of the $z_i$'s and a cutoff value.  If the statistic is above the cutoff we choose the specified distribution; if not, we choose the other.  While these strategies give lower bounds immediately, in Section \ref{upperboundsec} we will analyze the optimal strategy in each regime and show that these simple strategies are in fact asymptotically optimal and thus give the correct asymptotic total variation distance.

The benefit of these strategies is that they are simpler and more descriptive than the optimal strategy of comparing $PL(Z)$ to $ST(Z)$: they tell us what feature of the planted distribution takes the longest to `forget'.

\subsubsection{Flat Regime}
\label{flatlower}
We first describe the strategy in the first special case, one ball planted in each bin, then extend this to the general Flat regime. With one ball planted in each bin, a natural statistic would be to choose the ST distribution if any bin were empty and choose the PL distribution otherwise.  This strategy would separate the two distributions up to $m \sim n \log n$, but for higher scalings of $m$ every bin has at least one ball whp under the ST distribution so the statistic fails to differentiate the distributions. A better statistic is the number of pairs of balls in the same bin
\begin{equation*}
 PAIRS(Z) = \sum_{i=1}^n \binom{z_i}{2}
\end{equation*}
To see why this separates the distributions, compare the first $n$ balls under each distribution.  Under the ST distribution, the expected number of pairs that end up in the same bin is $\binom{n}{2} \frac{1}{n}$, while under the PL distribution $0$ pairs are in the same bin.  The $j$-th ball, $n<j\leq m$, is placed randomly in both ST and PL distributions and adds the same $\frac{j-1}{n}$ to the expectation of $\sum_{i=1}^n \binom{z_i}{2}$ and so the difference in the expectations remains precisely $\frac{n-1}{2}$.  We can exploit this discrepancy to give a lower bound for the total variation distance.  If we write  
\begin{equation}
\label{qi}
 q_i = \frac{n}{m} \left (z_i - \frac{m}{n} \right)
\end{equation}
we get
\begin{equation*}
 PAIRS(Z) = \frac{m^2}{2n^2} \left ( \sum_i q_i^2 \right)  + \frac{m^2}{2n} - \frac{m}{2} 
\end{equation*}
It will be convenient later to use the scaled and shifted statistic $\sum_i q_i^2$ instead of $PAIRS(Z)$.  

For a general initial planting in the Flat regime we adjust the statistic by adding weights to the $q_i$'s.  The weights are the $a_i$'s themselves and naturally arise from the analysis in Section \ref{lemsec}.  The statistic is: 
\begin{equation*}
 F_A(Z)  = \sum_{i=1}^n a_i q_i^2
\end{equation*}

From calculations in Section \ref{calcsec}, 
\begin{align*}
 \E_{ST} F_A &= \frac{kn}{m} +o(1)\\
\E_{PL} F_A &= \frac{kn}{m} - \frac{k^2n}{m^2} +o(1)
\end{align*}
and
\begin{equation*}
 \var_{ST}(F_A) \sim \var_{PL} (F_A) \sim \frac{2k^2n}{m^2}
\end{equation*}
Under our scaling $m \sim c k n^{1/2}$, and we set $\mu_F = \frac{kn}{m} - \frac{k^2n}{2m^2}$ to be the average of the two means.  Our strategy for distinguishing the two distributions is to choose the ST distribution if $F_A(Z) \ge \mu_F$ and choose the PL distribution if $F_A(Z) < \mu_F$.   We show in Section \ref{normalsec} that $F_A$ is asymptotically normal under each distribution.  Since the distance of $\mu_F$ from each mean is $\frac{\frac{k^2n}{2m^2} }{ \sqrt{ \frac{ 2k^2n}{m^2} } } = \frac{1}{2 \sqrt{2} c}$ standard deviations under the scaling $m \sim c n^{3/2}$, the asymptotic probability that $F_A(Z)$ is above or below $\mu_F$ in each case can be computed from the standard normal distribution function.  Thus the success probability for this strategy is:
\[ p = \frac{1}{2} \Phi \left( \frac{1}{2 \sqrt{2} c} \right) +   \frac{1}{2} \Phi \left( \frac{1}{2 \sqrt{2} c} \right) =   \Phi\left( \frac{1}{2 \sqrt{2} c} \right) \]
and so $ TV(m,n) \ge 2\Phi \left( \frac{1}{2 \sqrt{2} c}\right) -1 +o(1) $ from (\ref{tvlowerboundeq}), giving the lower bound for Theorem \ref{mainthm}.

\subsubsection{Hilly Regime}
\label{hillylower}
In the Hilly regime we define the statistic
\begin{equation*}
 H_A(Z) = \sum_{i=1}^n a_i q_i
\end{equation*}
with $q_i$ as in (\ref{qi}). If all $k$ balls are planted in bin $1$, $H_A(Z) = \frac{kn}{m} z_1 - k$, the number of balls in bin $1$ scaled and centered, so $H_A$ is equivalent in this case to the natural statistic mentioned in the first paragraph of this section. Calculations in Section \ref{calcsec} give
\begin{equation*}
 \E_{ST} ( H_A) = 0
\end{equation*}
 \begin{equation*}
 \E_{PL} (H_A) = \frac{ V n^2}{m} 
\end{equation*}
\begin{equation*}
 \var_{ST} (H_A) \sim \var_{PL} (H_A) \sim \frac{ V n^2}{m} 
\end{equation*}
Notice that the difference in expectations is completely accounted for by the first $k$ balls: if we let $\frac{n}{m}(a_i-1)$ be the contribution to $H_A$ of a ball that lands in bin $i$, then a randomly thrown ball contributes $0$ to the expectation, while the $k$ balls planted according to the initial planting $A$ contribute $\frac{n}{m} \sum_i a_i(a_i -1) = \frac{Vn^2}{m}$.  Also this choice of coefficients for a linear combination of the $q_i$'s maximizes the difference between the two means if we normalize by fixing the sum of the squares of the coefficients.

 We set $\mu_H = \frac{V n^2}{2m}$ to be the average of the two means, and our strategy is to choose the ST distribution if $H_A \ge \mu_H$ and the PL distribution otherwise.  Here the distance of $\mu_H$ from each mean is $\frac{1}{2 \sqrt c}$ standard deviations and $H_A$ is asymptotically normal under each distribution (Section \ref{normalsec}).  So calculating as above, $ TV(m,n) \ge 2\Phi \left( \frac{1}{2 \sqrt{c} }\right) -1 +o(1) $.

\subsubsection{Intermediate Regime}
\label{intlower}
In the Intermediate regime our statistic $I_A(Z)$ is a mixture of the two statistics in the previous regimes:
\begin{equation*}
 I_A(Z) = \sum_i a_i q_i - \frac{1}{2} \sum_ i a_i q_i^2
\end{equation*}
Under the scaling $m \sim cn^{3/2}$ with $V \sim \lam n^{-1/2}$, calculations give
\begin{equation*}
 \E_{ST} ( I_A) = - \frac{kn}{2m} + o(1)
\end{equation*}
 \begin{equation*}
 \E_{PL} (I_A) = - \frac{kn}{2m} + \frac{k^2n}{2m^2} + \frac{ V n^2}{m} + o(1) 
\end{equation*}
\begin{equation*}
 \var_{ST} (I_A) \sim \var_{PL} (I_A) \sim  \frac{V n^2}{m} + \frac{k^2n}{2m^2}
\end{equation*}
The average of the means is $\mu_I = - \frac{kn}{2m} + \frac{k^2 n}{4m^2} + \frac{V n^2}{2m}$ and the distance from each mean is
\begin{equation*}
 \frac{1}{2} \sqrt{\frac{\lam}{c} +\frac{1}{2c^2} } \, \,  \text{ standard deviations }
\end{equation*}
So we have
\begin{equation*}
 TV(m,n) \ge  2\Phi \left( \frac{1}{2} \sqrt{\frac{\lam}{c} +\frac{1}{2c^2} }\right) -1 +o(1)
\end{equation*}

\section{Upper Bounds}
\label{upperboundsec}

The statistics and strategies from the previous section give one way to distinguish between the two distributions, but in principle there could be better methods for distinguishing them.  In this section we show that the above strategies are asymptotically as good as the optimal strategy which consists of choosing the larger of $ST(Z)$ and $PL(Z)$.  We state a Lemma:

\begin{lem}
 \label{mainlem}
Let $V = \frac{ \sum a_i^2}{n}  - \frac{k^2}{n^2}$ as in (\ref{Vdef}). Then 
 \begin{enumerate}
\item \textbf{Flat Regime}: Let $m\sim c k n^{1/2}$ and assume $V = o\left ( \frac{k}{n^{3/2}} \right)$. Set $\mu_F = \frac{kn}{m} - \frac{k^2n}{2m^2}$ as in Section \ref{flatlower}.  For any fixed $\eps >0$,  with ST and PL probability $1-o(1)$,
\begin{enumerate}
\item $F_A(Z) \ge \mu_F + \eps  \Rightarrow ST(Z) > PL(Z)$
\item $F_A(Z) \le \mu_F - \eps  \Rightarrow PL(Z) > ST(Z) $
\end{enumerate}

\item \textbf{Hilly Regime}: Let $m\sim c V n^2$ and assume $V = \omega \left ( \frac{k}{n^{3/2}} \right)$.  Set $\mu_H = \frac{Vn^2}{2m}$ as Section \ref{hillylower}.  For any fixed $\eps >0$ the following hold with ST and PL probability $1-o(1)$:
\begin{enumerate}
\item $ H_A(Z) \ge \mu_H + \eps  \Rightarrow PL(Z) > ST(Z)$
\item $H_A(Z) \le \mu_H - \eps  \Rightarrow ST(Z) > PL(Z) $
\end{enumerate}

\item \textbf{Intermediate Regime}: Let $m \sim c k n^{1/2}$ and assume $V \sim \lam \frac{k}{n^{3/2}}$.  Set $\mu_I =- \frac{kn}{2m} + \frac{k^2 n}{4m^2} + \frac{V n^2}{2m}$ as in Section \ref{intlower}.  Then for any fixed $\eps >0$, with ST and PL probability $1-o(1)$,
\begin{enumerate}
\item $I_A(Z) \ge \mu_I + \eps  \Rightarrow PL(Z) > ST(Z) $
\item $I_A(Z) \le \mu_I - \eps \Rightarrow ST(Z) > PL(Z) $
\end{enumerate}

\end{enumerate}
\end{lem}

We will prove Lemma \ref{mainlem} in Section \ref{lemsec}. The proof of Theorem \ref{mainthm} given Lemma \ref{mainlem} is very similar in all three cases.  Here we prove it for the Hilly regime and omit the Flat and Intermediate cases.
\begin{proof}[Proof of Theorem \ref{mainthm}, Hilly case] Fix some $\eps > 0$.  We partition the set of configurations as follows:
\begin{align*}
\Omega_1 &= \{ Z: ST(Z) > PL(Z) \text { and } H_A(Z) \le \mu_H - \eps  \} \\
\Omega_2 &= \{ Z: PL(Z) > ST(Z) \text { and } H_A(Z) \ge \mu_H + \eps  \} \\
 \Omega_3 &= \{ Z: ST(Z) > PL(Z) \text { and } H_A(Z) \ge \mu_H + \eps \} \\
 \Omega_4 &= \{ Z: PL(Z) > ST(Z) \text { and } H_A(Z) \le \mu_H - \eps \} \\
 \Omega_5 &= \{ Z: H_A(Z) \in ( \mu_H -\eps , \mu_H + \eps) \}
\end{align*}
For $n$ large, $ST(\Omega_5) \le \frac{2}{\sqrt c} \eps$ and $PL(\Omega_5) \le \frac{2}{\sqrt c} \eps$ since $H_A$ is asymptotically normal and the interval $ ( \mu_H -\eps, \mu_H + \eps)$ has width $\frac{2}{\sqrt c} \eps$ standard deviations under either the ST or PL distribution.  Now recall the definition of total variation distance from equation (\ref{tvdef}). This is equivalent to 
\begin{equation*}
 ||ST -PL ||_{TV} = \sum_{Z: ST(Z) > PL(Z)} ST(Z) - PL(Z) 
\end{equation*}

Define the similar quantity
\begin{equation*}
 TV^\prime(m,n) = \sum_{Z: H_A(Z) < \mu_H } ST(Z) - PL(Z)
\end{equation*}
which is what we would get as an estimate for the total variation distance if we used our strategy from the lower bound.  We now show that $TV(m,n) =  TV^\prime(m,n) +o(1)$.

\begin{equation*}
 TV(m,n) \le  \sum_{\Omega_1} ST(Z) - PL(Z) + \sum_{\Omega_3} ST(Z) + \sum_{\Omega_5} ST(Z)
\end{equation*}
and 
\begin{equation*}
 TV^\prime(m,n) \ge \sum _{\Omega_1} ST(Z) - PL(Z) - \sum_{\Omega_4} PL(Z) - \sum_{\Omega_5} PL(Z)
\end{equation*}
So
\begin{equation*}
 TV(m,n) \le TV^\prime(m,n)+ ST(\Omega_3) +  PL(\Omega_4)  +  \sum_{\Omega_5} ST(Z) + PL(Z)
\end{equation*}
$ST(\Omega_3)$ is $o(1) $ by Lemma \ref{mainlem} part 2a.  $PL(\Omega_4)$ is $o(1) $ by Lemma \ref{mainlem} part 2b.   The sum over $\Omega_5$ is $\le \frac{4}{\sqrt c} \eps$.  So $TV(m,n) \le  TV^\prime(m,n) + \frac{4}{\sqrt c} \eps +o(1)$.  Since $\eps$ is arbitrary, $TV(m,n) \le TV^\prime(m,n) +o(1)$.  The other side of the inequality is similar and so $ TV(m,n) =  TV^\prime(m,n) + o(1)$.   \end{proof}

\subsection{Proof of Upper Bound Lemma}
\label{lemsec}
The idea of Lemma \ref{mainlem} is that with high probability over the choice of configuration, each of our simple lower bound strategies based on a threshold statistic is in fact equivalent to the optimal strategy.  The optimal strategy is to compute $\frac{PL(Z)}{ST(Z)}$ and pick PL if the ratio is $\ge 1$.  To prove the Lemma we compute the logarithm of the ratio, expand the terms and show that in each particular regime all the terms are concentrated except for one term which corresponds to the respective statistic.  Then the sign of the logarithm is determined by whether the statistic is above or below its threshold.  

\begin{proof} The exact formulae for $ST(Z)$ and $PL(Z)$ are 
\begin{align*}
 ST(Z) &= \frac{ 1}{n^m} \binom{ m}{z_1 \dots z_n} \\
 PL(Z) &= \frac{1}{n^{m-k}} \binom { m-k  } {(z_1 - a_1) \cdots (z_n - a_n) } 
\end{align*}
and so
\begin{equation}
\label{ratioeq}
 \frac{ PL(Z)}{ST(Z)} = \frac{n^k (z_1)_{a_1} \cdots (z_n)_{a_n}  } { (m)_k  }
\end{equation}
 We write
\begin{align*}
 \prod_i (z_i)_{a_i} &= E_1 \cdot \prod_i z_i^{a_i} \\
(m)_k = \frac{m^k}{E_2} 
\end{align*}
and get
\begin{equation*}
 \frac{ PL(Z)}{ST(Z)} \sim  E_1 E_2  \prod_i (1+q_i)^{a_i}
\end{equation*}
where $q_i = \frac{n}{m}(z_i - \frac{m}{n})$. In Section \ref{errorsec} we show that  
\begin{equation*}
  \ln (E_1 E_2)  = -\frac{V n^2}{2m} + \frac{kn}{2m}+ \frac{5 kn^2}{12 m^2} - \frac{k^2 n }{4m^2} +o(1)
\end{equation*}
 whp in each regime.  Taking the logarithm of the ratio, 
 \begin{equation*}
 \ln \frac{ PL(Z)}{ST(Z)} = -\frac{V n^2}{2m} + \frac{kn}{2m}+ \frac{5 kn^2}{12 m^2} - \frac{k^2 n }{4m^2} + \sum_{i=1}^n a_i \log (1+ q_i) + o(1) 
\end{equation*}
Chernoff bounds give that the $q_i$'s are uniformly $o(1)$ whp in all three regimes, so we can use the Taylor series for the logarithm:
\begin{align*}
\ln \frac{ PL(Z)}{ST(Z)} &= -\frac{V n^2}{2m} + \frac{kn}{2m}+ \frac{5 kn^2}{12 m^2} - \frac{k^2 n }{4m^2}\\
 &+\sum_i a_i q_i - \frac{1}{2} \sum _i a_i q_i^2 + \frac{1}{3}\sum_i a_i q_i^3  -\frac{1}{4}  \sum_i a_i q_i^4  + o(1)  
\end{align*}
where the higher order terms are $o(1)$ from the Chernoff bound.  Next we see from our calculations in Section \ref{calcsec} that the variances of $\sum a_i q_i^3$ and $\sum a_i q_i^4$ are $o(1)$ in all regimes, and so they are concentrated around their means.  This gives whp:
\begin{equation}
\label{logratioeq}
\ln \frac{ PL(Z)}{ST(Z)} = -\frac{V n^2}{2m} + \frac{kn}{2m} - \frac{k^2 n }{4m^2}+\sum_i a_i q_i - \frac{1}{2} \sum _i a_i q_i^2  + o(1)  
\end{equation}
 We now analyze this sum under the specifics of each regime.

\subsubsection{Flat Regime}

Here $m \sim c k n^{1/2}$ and $V = o\left ( \frac{k}{n^{3/2}}\right)$.  Here the mean and variance of $\sum _i a_i q_i$ are $o(1)$ for both the ST and PL distributions, and $\frac{V n^2}{2m} = o(1)$, so whp we have:
\begin{equation}
 \ln \frac{ PL(Z)}{ST(Z)}  = \frac{kn}{2m} - \frac{k^2 n }{4m^2} - \frac{1}{2} \sum _i a_i q_i^2  + o(1)
\end{equation}

 But $\sum _i a_i q_i^2$ is precisely our statistic $F_A(Z)$.  By assumption in Lemma \ref{mainlem} part 1a), $F_A(Z) \ge \mu_F + \eps = \frac{kn}{m} - \frac{k^2n}{2m^2} + \eps$, and $ - \frac{1}{2} F_A(Z) \le - \frac{kn}{2m} + \frac{k^2n}{4m^2} + - \frac{\eps}{2}$ so
\begin{equation}
 \ln \frac{ PL(Z)}{ST(Z)} \le -\frac{\eps}{2} + o(1)
\end{equation}
and so whp, $ST(Z) > PL(Z)$.  Proving 1b) is similar: the assumption says that $F_A(Z) \le \frac{kn}{m} - \frac{k^2 n}{2m^2} - \eps$, which gives $\ln \frac{ PL(Z)}{ST(Z)} \ge \frac{\eps}{2} + o(1)$, and so $PL(Z) > ST(Z)$ whp.

\subsubsection{Hilly Regime}

In the Hilly regime we have $V = \omega(k n^{-3/2})$ and $m \sim c V n^2$.  In this regime $\frac{k^2 n}{m^2} = o(1)$ and $\sum a_i q_i^2 $ is concentrated around its mean, giving simply:
\begin{equation*}
 \ln \frac{ PL(Z)}{ST(Z)}  = - \frac{Vn^2}{2m}  + \sum_i a_i q_i + o(1)
\end{equation*}
whp. $\sum_i a_i q_i$ is $H_A(Z)$.  By assumption in part 2a) of Lemma \ref{mainlem}, $H_A \ge \frac{Vn^2}{2m} + \eps$, so $\ln \frac{ PL(Z)}{ST(Z)} \ge \eps +o(1)$ and $ST(Z) > PL(Z)$ whp.  Similarly under the assumptions of 2b) we get $PL(Z) > ST(Z)$ whp.

\subsubsection{Intermediate Regime}

The intermediate regime is similar:  $\sum_i a_i q_i - \frac{1}{2} \sum _i a_i q_i^2$ is the statistic $I_A(Z)$, and the conditions on $I_A$ imply the result. \end{proof}

\section{Asymptotic Normality}
\label{normalsec}

\subsection{Hilly Regime}

In the Hilly regime, our statistic $H_A(Z) = \sum_{i=1}^n a_i q_i$ under the ST distribution can be rewritten as $ H_A(Z) = \sum_{j=1}^m Y_j$ where the $Y_j$'s are i.i.d., one for each ball, with $Y_j = \frac{n}{m}(a_i - \frac{k}{n})$ if ball $j$ is in bin $i$.  To apply the Lindeberg-Feller Central Limit Theorem (see, for example, \cite{durrett2010probability}), we check that the $Y_i$'s are bounded.
$Y_j \ge - \frac{k}{m} = o(1)$, and
\begin{align*}
 Y_j \le \frac{nM}{m}
\end{align*}
where $M = \max_i a_i$. $\sum a_i^2 = Vn + \frac{k^2}{n^2}$, so $M \le \sqrt {\frac{k^2}{n^2} + Vn}$ and
 \begin{align*}
 Y_j &\le \frac{\sqrt{ k^2 + Vn^3}}{cVn^2} \\
&\le 2\frac{k}{ c n^2 V} + 2\sqrt{ \frac{1}{c^2 n V}} = o(1)
\end{align*} 
since we are in the Hilly regime and $k>> \sqrt n$. The Lindeberg-Feller CLT then implies that $\frac{ H_A}{\sqrt{ \var(H_A)}} \Rightarrow N(0,1)$.  

Under the PL distribution we write $H_A(Z) = \frac{n}{m} \sum a_i^2 + \sum_{j=1}^{m-n} Y_j$ where the $Y_j$'s are as above and correspond to the $m-n$ randomly placed balls.  Again Lindeberg-Feller gives $\frac{ H_A - \E_{PL} H_A}{\sqrt{ \var(H_A)}} \Rightarrow N(0,1)$.  

\subsection{Flat and Intermediate Regime}

In the Flat and Intermediate regimes there are not such simple representations of $F_A$ or $I_A$ as sums of independent random variables, but both statistics fit into a general framework of statistics of `occupancy scores' that have been proved to have normal limits is a series of papers (\cite{morris1966admissible}, \cite{morris1975central}, \cite{quine1982berry}, \cite{quine1984normal}).  The following Theorem from \cite{morris1966admissible} suffices for our cases:

\begin{thm}
 \label{normalthm}
Let $z_1, \dots z_n$ be a multinomial vector with parameters $(p_1, \dots p_n)$ and $\sum_i z_i = m$. Let $f_1, \dots f_n$ be degree 2 polynomials and $S_n = \sum_{i=1}^n f_i(z_i)$.  Suppose:
\begin{enumerate}
 \item  $\max_{1 \le i \le n} p_i = o(1)$
\item $\min_{1 \le i \le n} m p_i $ is bounded away from $0$ as $n \to \infty$
\item $\max_{1 \le i \le n} \var(f_i(z_i)) / \sum_1^n \var(f_j(z_j)) = o(1)$
\end{enumerate}

Then $\frac{S_k - \E S_k}{\sqrt{\var(S_k)}} \Rightarrow N(0,1)$.
\end{thm}

In our case $p_i = \frac{1}{n}$ for all $i$, and $\frac{m}{n}$ is bounded away from $0$ under our scalings.  Under the ST distribution we set $f_i(z_i) = \frac{n^2}{m^2} a_i \left( z_i - \frac{m}{n} \right)$, a polynomial of degree 2.  The $z_i$'s have the same distribution and the $f_i$'s are the same except for the factor $a_i$, so $ \frac{\var(f_i(z_i)) }{ \sum_1^n \var(f_i(z_i))} = \frac{a_i^2}{\sum a_j^2}$.  In the Flat and Intermediate regime $\max_i a_i \le \frac{k}{n} + \sqrt{Vn}$, and $\sum a_i^2 \ge \frac{k^2}{n}$, so $\frac{a_i^2}{\sum a_j^2} = o(1)$ and Theorem \ref{normalthm} applies.

For the PL distribution we let $y_i$ be the number of the $m-k$ random balls that end up in bin $i$. Here the $y_1, \dots y_n$ are a multinomial vector with $p_i = \frac{1}{n}$ and $\sum y_i = m -k$.  We let $f_i(z_i) = \frac{n^2}{m^2} a_i \left ( a_ i + y_i - \frac{m}{n} \right )^2$.   As above, $\frac{\var(f_i(z_i)) }{ \sum_1^n \var(f_i(z_i))} = o(1)$ as long as $\var \left (  ( a_ i + y_i - \frac{m}{n}  )^2 \right) \sim \var \left (  ( a_ j + y_j - \frac{m}{n}  )^2 \right)$ for all $i, j$.  This holds since $\max a_k = o\left(\frac{m}{n} \right)$ and both variances become $\var \left (y_i^2 - 2 \frac{m}{n} y_i (1+ o(1)) \right)$.

\section{Error Term}
\label{errorsec}
We will need the following proposition about the $a_i$'s:
\begin{prop}
\label{aiprop}
 In all three regimes,
\begin{enumerate}
 \item $\sum a_i^2 = \frac{k^2}{n} + nV$
\item $\sum a_i^3 = \frac{k^3}{n^2} + o\left(\frac{m^2}{n^2} \right)$
\item $\sum a_i^4 = \frac{k^4}{n^3} + o\left(\frac{m^3}{n^3}\right )$
\end{enumerate}
\end{prop}

\begin{proof}
 1) is the definition of $V$.  For 2), we write $a_i = \frac{k}{n} + \Delta_i$.  Then we have
\begin{align*}
 \sum_{i=1}^n a_i^2 &= \sum \frac{k^2}{n^2} + 2 \frac{k}{n} \Delta_i + \Delta_i^2 \\
&= \frac{k^2}{n} + \sum \Delta_i^2
\end{align*}
where we use the fact that $\sum \Delta_i =0$, so $\sum \Delta_i^2 = Vn$ from the definition of $V$.  Now we write
\begin{align*}
 \sum _{i=1}^n a_i^3 &= \sum_{i=1}^n \frac{k^3}{n^3}  + 3 \frac{k^2}{n^2} \Delta_i + 3 \frac{k}{n} \Delta_i^2 + \Delta_i^3 \\
&= \frac{k^3}{n^2} + 3 k V + \sum_{i=1}^n \Delta_i^3
\end{align*}
It is straightforward to check that $kV = o \left( \frac{m^2}{n^2} \right )$ in all three regimes.  We bound 
\begin{equation*}
 \sum_{i=1}^n \Delta_i^3 \le \left( \sum_{i=1}^n a_i^2 \right )^{3/2} = (Vn)^{3/2}
\end{equation*}
and similarly we can check that $(Vn)^{3/2} = o \left( \frac{m^2}{n^2} \right )$ in all three regimes.

For 3), we write
\begin{align*}
 \sum _{i=1}^n a_i^4 &= \sum_{i=1}^n \frac{k^4}{n^4}  + 4 \frac{k^3}{n^3} \Delta_i + 6\frac{k^2}{n^2} \Delta_i^2 + 4 \frac{k}{n} \Delta_i^3 + \Delta_i^4 \\
&= \frac{k^4}{n^3} + 6 \frac{k^2}{n} V + 4 \frac{k}{n} \sum_{i=1}^n \Delta_i^3 + \sum_{i=1}^n \Delta_i^4
\end{align*}
Again we can check that $\frac{k^2 V}{n} = o \left( \frac{m^3}{n^3} \right)$ in all three regimes, and we bound the last two terms by $\frac{k}{n} (Vn)^{3/2}$ and $(Vn)^2$ respectively, both of which $= o \left( \frac{m^3}{n^3} \right)$ in all three regimes. 
\end{proof}

Now we prove the main Lemma of this section:

\begin{lem}
\label{kElem}
 Let 
\begin{equation*}
 E_1 = \frac{\prod_{i=1}^n (z_i)_{a_i} }{\prod_{i=1}^n z_i^{a_i} } 
\end{equation*}
and
\begin{equation*}
 E_2 = \frac{m^k}{(m)_k}
\end{equation*}

Then whp in all regimes,
\begin{equation*}
  \ln (E_1 E_2) = -\frac{V n^2}{2m} + \frac{kn}{2m}+ \frac{5 kn^2}{12 m^2} - \frac{k^2 n }{4m^2}
\end{equation*}
\end{lem}

First we use a standard asymptotic approximation of $(m)_k$:
\begin{equation*}
\ln  E_1 = \ln \frac{m^k}{(m)_k} =\frac{k^2}{2m} + \frac{k^3}{6m^2} + \dots 
\end{equation*}

Next we compute the asymptotics of $E_2$.
\begin{align*}
  E_2 &= \prod_{i=1}^n \frac{ (z_i)_{a_i} }{ z_i^{a_i}} \\
&= \prod_{i=1}^n 1 \cdot \left(1- \frac{1}{z_i}\right ) \cdots \left (1 - \frac{a_i-1}{z_i} \right) 
\end{align*}
and so
\begin{align*}
 \ln E_2 &= \sum_{i=1}^n \sum_{j=0}^{a_i-1} \ln \left (1 - \frac{j}{z_i} \right ) \\
&= -\sum_{i=1}^n \sum_{j=0}^{a_i-1}  \frac{j}{z_i} + \frac{j^2}{2z_i^2} + \frac{j^3}{3z_i^3} + \dots \\
&= - \sum_{i=1}^n \left [ \frac{1}{z_i} \sum_{j=0}^{a_i-1} j + \frac{1}{2z_i^2} \sum_{j=0}^{a_i-1} j^2 +  \frac{1}{3z_i^3} \sum_{j=0}^{a_i-1} j^3 +\dots \right ] \\
&= - \sum_{i=1}^n \left [ \frac{1}{z_i} \frac{a_i^2 -a_i}{2} + \frac{1}{2z_i^2} \frac{a_i (a_i-1)(2 a_i-1) } {6 } + \frac{1}{3z_i^3} \frac{a_i^2 (a_i-1)^2 } {4 }+ \dots \right ] 
\end{align*}

Using long division we see that
\begin{equation*}
 \frac{1}{z_i} = \frac{n}{m} \left (1 -q_i + q_i^2 - q_i^3 + \dots  \right ) 
\end{equation*}
and a Chernoff bound shows that whp
\begin{equation*}
 z_i - \frac{m}{n} = O \left(\sqrt {\log n} \sqrt{ \frac{m}{n}}   \right )
\end{equation*}
for all $i$. We expand the expression for $\ln E_2$ term by term:
\begin{align*}
 - \sum_{i=1}^n \frac{1}{z_i} \frac{a_i^2 - a_i}{2} &= -\frac{n}{2m} \left( Vn + \frac{k^2}{n} -   \sum_{i=1}^n( a_i^2 -a_i)(q_i-q_i^2 + \dots  )  \right) \\
&= - \frac{Vn^2}{2m} - \frac{k^2}{2m} - \frac{n}{2m} \sum_{i=1}^n a_i^2 q_i - a_i q_i + a_i^2 q_i^2 - a_i q_i^2 + \dots
\end{align*}
Using Proposition \ref{aiprop} and calculations from Section \ref{calcsec}, we see that whp, $\frac{n}{2m} \sum a_i^2 q_i^2 = \frac{k^2n}{2m^2}+o(1)$, $\frac{n}{2m}\sum a_i q_i^2 = \frac{kn^2}{2m^2}$, and the other terms in the sum are all $o(1)$.  This gives:

\begin{equation*}
 - \sum_{i=1}^n \frac{1}{z_i} \frac{a_i^2 - a_i}{2} =- \frac{Vn^2}{2m} - \frac{k^2}{2m} - \frac{k^2n}{2m^2} + \frac{kn^2}{2m^2} + o(1)
\end{equation*}
The next term is:
\begin{align*}
- \sum_{i=1}^n \frac{1}{2z_i^2} \frac{a_i (a_i-1)(2 a_i-1) } {6 } &=  - \frac{n^2}{12m^2} \sum_{i=1}^n (2a_i^3 - 3 a_i^2 + a_i) (1 - 2q_i+ 3q_i^2 - \dots)
\end{align*}
Again using Proposition \ref{aiprop}, we calculate that whp
\begin{align*}
 - \sum_{i=1}^n \frac{1}{2z_i^2} \frac{a_i (a_i-1)(2 a_i-1) } {6 } &= - \frac{k^3}{6 m^2} + \frac{ k^2n}{4m^2} - \frac{kn^2}{12 m^2} + o(1)
\end{align*}
For the remaining terms, only the leading terms will be $\Theta(1)$, and that only for large values of $k$.  Those leading terms are precisely the negative of the remaining terms in the asymptotics $\ln \frac{m^k}{(m)_k} =\frac{k^2}{2m} + \frac{k^3}{6m^2} + \dots$.  Putting this together gives the asymptotics of $\ln (E_1 E_2)$ and Lemma \ref{kElem}.

\section{Calculations}
\label{calcsec}
In this section we calculate the mean and variance of $\sum a_i q_i$, $\sum a_i q_i^2$, $\sum a_i q_i^3$, and $\sum a_i q_i^4$.   The results:

\begin{prop}
\label{kaisprop}
In each of our three regimes the following hold:
\begin{itemize}
\item $\E_{ST} \sum a_i q_i = 0$, $\var_{ST} ( \sum a_i q_i) = \frac{V n^2}{m}$
\item $\E_{PL} \sum a_i q_i = \frac{V n^2}{m}$, $\var_{PL} ( \sum a_i q_i) = \frac{V n^2}{m} +o(1)$

\item $\E_{ST} \sum a_i q_i^2 = \frac{kn}{m} +o(1) $, $\var_{ST} ( \sum a_i q_i^2) = \frac{2 k^2n}{m^2} +o(1)$
\item $\E_{PL} \sum a_i q_i^2 = \frac{kn}{m} - \frac{k^2n}{m^2}   + o(1) $, $\var_{PL} ( \sum a_i q_i^2) = \frac{2 k^2 n}{m^2} +o(1)$

\item $\E_{ST} \sum a_i q_i^3 = \frac{kn^2}{m^2} +o(1)$, $\var_{ST} ( \sum a_i q_i^3) = o(1)$
\item $\E_{PL} \sum a_i q_i^3 = \frac{kn^2}{m^2} +o(1)$, $\var_{PL} ( \sum a_i q_i^3)  =o(1)$

\item $\E_{ST} \sum a_i q_i^4 = \frac{3kn^2}{m^2} +o(1)$, $\var_{ST} ( \sum a_i q_i^4) = o(1)$
\item $\E_{PL} \sum a_i q_i^4 = \frac{3kn^2}{m^2} +o(1)$, $\var_{PL} ( \sum a_i q_i^4) = o(1)$
\end{itemize}

\end{prop}

We show the calculations for $\sum a_i q_i$ and $\sum a_i q_i^2$.  The calculations for $\sum a_i q_i^3$ and $\sum a_i q_i^4$ are somewhat tedious and unenlightening and so are omitted. 

\subsection{$\sum a_i q_i$}
 For the ST distribution: We write $\sum _i a_i q_i = \sum_{j=1}^m Y_j$ where the $Y_j$'s are i.i.d. and $Y_j = \frac{n}{m}(a_i - \frac{k}{n})$ if ball $j$ is in bin $i$.  From this we see $\E_{ST} \sum a_i q_i =0$ and $\var_{ST}( \sum a_i q_i ) = m \cdot \var (Y_j) = m \frac{V n^2}{m^2} = \frac {V n^2}{m}$.

For the PL-distribution, we write $\sum _i a_i q_i = \frac{n}{m} ( \sum a_i(a_i - \frac{k}{n} )+ \sum_{j=1}^{m-k} Y_j = \frac{Vn^2}{m} + \sum_{j=1}^{m-k} Y_j$ where the $Q_j$'s are the same as above but indexed over the $m-k$ randomly placed balls.  This gives $\E_{PL}(\sum a_i q_i) = \frac{Vn^2}{m}$ and $\var_{PL}(\sum a_i q_i) = (m-k) \cdot \var(Y_j) \sim \frac{Vn^2}{m}$.

\subsection{$\sum a_i q_i^2$}

for the ST distribution we write
\begin{equation*}
  \sum _i a_i q_i^2 = k + \sum_{j <k} P_{j,r} + \sum _{l=1}^m Y_l
\end{equation*}
where $P_{j,r} = \frac{2 n^2}{m^2} a_i$ if the pair of balls $j,r$ are both in bin $i$, $P_{j,r} = 0$ if $j, k$ are in different bins, and $Y_l = \left( \frac{n^2}{m^2} - \frac{2n}{m}  \right) a_i $ if ball $l$ is in bin $i$. So
\begin{equation*}
 \E_{ST} \sum _i a_i q_i^2 = k + \binom{m}{2}  \E P_{j,r} + m  \E Y_l
\end{equation*}
We calculate that $ \E_{ST} P_{j,r} =  \frac{2k}{m^2} $ and $\E _{ST} Y_l = \frac{kn}{m^2} -  \frac{2k}{m}$ giving
\begin{equation*}
 \E_{ST}  \sum _i a_i q_i^2 = \frac{kn}{m}  - \frac{k}{m} = \frac{kn}{m} + o(1) 
\end{equation*}
And the variance is:
\begin{align*}
 \var_{ST} ( \sum a_i q_i^2 ) = &\binom{m}{2}  \var (P_{j,r}) \\
 &+ m \cdot  \var( Y_l)  + (m)_3  \cov( P_{j,r}, P_{r,l} )+ 4\binom{m}{2}  \cov (P_{j,l}, Y_l)
\end{align*}
since the $P_{j,r}$'s and $Y_l$'s that do not share an index are independent.  We calculate
\begin{align*}
 \var (P_{j,r}) &= \sum_i a_i^2 \frac{4 n^2}{m^4} - \frac{4 k^2}{m^4} \\
&=  \frac{4n k^2}{m^4} +  \frac{4 V n^3}{m^4} - \frac{4 k^2}{m^4} \\
 \var(Y_l) &= \sum_i a_i^2  \left ( \frac{n^3}{m^4} -  \frac{4n^2}{m^3} + \frac{4n}{m^2}   \right )  - \frac{k^2n^2}{m^4} +  \frac{4k^2n}{m^3} -  \frac{4k^2}{m^2}\\
 &= \frac{V n^4}{m^4} -  \frac{4 V n^3}{m^3} +  \frac{4V n^2}{m^2} \\
 \cov( P_{j,r}, P_{r,l} ) &= \sum_ i  \frac{ 4 n}{m^4} a_i^2 -   \frac{4k^2}{m^4} = 4\left (\frac{k^2}{n} + Vn \right) \frac{n}{m^4} -   \frac{4k^2}{m^4} \\
 &=  \frac{4 V n^2 }{m^4}\\
 \cov (P_{j,l}, Y_l) &=   \sum_i   \left ( \frac{2 n^2}{m^4} -  \frac{ 4 n}{m^3} \right ) a_i ^2 -  2 \frac{k^2 n}{m^4} +  \frac{4k^2}{m^3} \\
&=   \frac{ 2 V n^3}{m^4} - \frac{4V n^2}{m^3}  
\end{align*}
which all together gives
\begin{equation*}
 \var_{ST} \left ( \sum a_i q_i^2 \right) =   \frac{2 k^2n}{m^2} + o(1)  
\end{equation*}
in all our regimes as long as $k >> \sqrt n$.

\vf

And for the PL distribution:
\begin{equation*}
 \sum_i a_i q_i^2 = \frac{n^2}{m^2} \sum _i a_i \left (a_i + r_i - \frac{m}{n} \right )^2 
\end{equation*}
where $r_i$ is the number of randomly thrown balls in bin $i$, so $z_i = a_i + r_i$.  Expanding we find the constant part is:
\begin{equation*}
 C= \frac{n^2}{m^2} \sum_i a_i \left ( a_i^2  + \frac{m^2}{n^2}     - 2 a_i \frac{m}{n} \right )  =  \frac{k^3}{m^2}  +k -  \frac{2k^2}{m} -  \frac{2V n^2}{m} +o(1)
\end{equation*}
since in all regimes $ \frac{ n^2 \sum a_i^3}{m^2} = \frac{k^3}{m^2} + o(1)$ (see Proposition \ref{kaisprop}).  The random part is:
\begin{equation*}
 \frac{n^2}{m^2} \sum_i a_i r_i^2 + \frac{n^2}{m^2} \sum _i \left( 2 a_i^2 - 2 a_i \frac{m}{n}   \right ) r_i
\end{equation*}
If we let $R_i$ be the number of pairs of random balls in bin $i$, then $r_i ^2 = 2 R_i + r_i$ and the random part becomes:
\begin{equation*}
 \frac{ 2n^2}{m^2} \sum_i a_i R_i + \frac{n^2}{m^2} \sum_i \left( 2 a_i^2 - 2 a_i \frac{m}{n}  +a_i  \right ) r_i
\end{equation*}
So we define $P_{j,r}$ to be $\frac{2n^2}{m^2} a_i$ when the pair of random balls $j, r$ are both in bin $i$, and $Y_l =  \frac{n^2}{m^2} \left( 2 a_i^2 - 2 a_i \frac{m}{n}  +a_i  \right )$ when ball $l$ is in bin $i$, and then 
\begin{equation*}
 \sum a_i q_i^2 = C + \sum_{j < r} P_{j,r} + \sum_{l=1}^{m-k} Y_l
\end{equation*}
We have:
\begin{align*}
 \E Y_l &=   \frac{2k^2}{m^2} + \frac{kn}{m^2} + \frac{2 V n^2}{m^2} -  \frac{2k}{m} \\ 
 \E P_{j,r} &=  \frac{2k}{m^2}
\end{align*}
and so
\begin{align*}
 \E_{PL} \sum _i a_i q_i ^2 = &\frac{k^3}{m^2}  +k -  \frac{2k^2}{m} - \frac{2Vn^2}{m} \\
&+ (m-k)   \left(   \frac{2k^2}{m^2} + \frac{kn}{m^2} + \frac{2 V n^2}{m^2} -  \frac{2k}{m}     \right ) + \binom{m-k}{2}  \frac{2k}{m^2}       +o(1) 
\end{align*}
\begin{equation*}
 = \frac{kn}{m} - \frac{k^2n}{m^2}   + o(1) 
\end{equation*}
since $\frac{Vn^2k}{m^2}$ and $\frac{k}{m} = o(1)$ in all regimes.

Variance:  The $P_{j,k}$'s are the same for the PL distribution as for the ST distribution and for the $Y_l$'s we calculate 
\begin{equation*}
 \var(Y_l) =    \frac{4V n^2}{m^2} + o\left (m^{-1} \right)
\end{equation*}
\begin{equation*}
 \cov (P_{j,l}, Y_l) = - \frac{4 V n^2}{m^3} + o \left (m^{-2} \right)   
\end{equation*}
which gives
\begin{equation*}
 \var_{PL} ( \sum a_i q_i^2 ) =  \frac{2 k^2n}{m^2} + o(1) 
\end{equation*}

\section{Concluding Remarks}
\label{concludesec}
We discuss briefly a natural modification of this problem suitable for further study.

\paragraph{Unlabeled Bins}

Suppose the bins are now indistinguishable.  It becomes more difficult now to distinguish the two distributions (since we could have ignored the labels in the previous case), so the correct scaling for $m$ will be no greater than in the original labeled case.  In the case of planting exactly one ball in each bin, the scaling is actually the same - as mentioned above, the statistic $\sum_i q_i^2$ is an asymptotically correct statistic in this regime, and since it does not depend on the labeling of the bins, we can use it just as is in the unlabeled case.  The Hilly regime however is much different.  There our distinguishing statistic, $\sum a_i q_i$ depended very much on the labeling of the bins.  In the case of $n$ balls planted in one bin and none in the rest, a natural guess of the right statistic would be the maximum number of balls in a bin.  This turns out to be correct and it comes with a different scaling for $m$, $m \sim \frac{n^3}{2 \log n} \left ( 1+ \frac{c}{\sqrt{ \log n}} \right )$ with $c \in (- \infty, \infty)$ a constant.  The distribution of the maximum number of balls in a bin, studied in \cite{raab1998balls}, is very useful in this case.  The general case seems to be related asymptotically to the problem of distinguishing between two scenarios: one in which we have $n$ independent $N(0,1)$ random variables and one in which we have a collection of normal random variables with means different than 0.  Such a problem was studied in a different context in \cite{donoho2004higher}.

\bibliographystyle{acm}	
\bibliography{ballsbins}

\end{document}